\documentclass[12pt]{article}
\usepackage[amsmath]{e-jc}

\usepackage{graphicx}
\usepackage{algorithm}
\usepackage{algorithmic}

\Copyright{The author. Released under the CC BY-ND license (International 4.0).}

\title{Ramsey Number Counterexample Checking and One Vertex Extension Linearly Bound by $s$ and $t$}

\author{Adam Lehavi\authornote{1}}

\authortext{1}{Department of Mathematics, University of Southern California, Los Angeles, United States
   (\email{alehavi@usc.edu}).}

\begin{document}

\maketitle

\begin{abstract}
The Ramsey number $R(s,t)$ is the smallest integer $n$ such that all graphs of size $n$ contain a clique of size $s$ or an independent set of size $t$. $\mathcal{R}(s,t,n)$ is the set of all counterexample graphs without this property for a given $n$. I prove that if a graph $G_{n+1}$ of size $n+1$ has $\max\{s,t\}+1$ subgraphs in $\mathcal{R}(s,t,n)$, then $G_{n+1}$ is in  $\mathcal{R}(s,t,n+1)$. Based on this, I introduce algorithms for one-vertex extension and counterexample checking with runtime linearly bound by $s$ and $t$. The utility of these algorithms is shown by verifying $\mathcal{R}(4,6,36)$ and $\mathcal{R}(5,5,43)$ are empty given current sets $\mathcal{R}(4,6,35)$ and $\mathcal{R}(5,5,42)$.
\end{abstract}

\section{Introduction}
The Ramsey number $R(s,t)$ is defined as the smallest positive integer $n$ such that all graphs of size $n$ contain a clique of size $s$ or an independent set of size $t$. The field of Ramsey theory explores finding exact values, upper bounds and lower bounds for given $s$ and $t$ values. The field began with the proof that $R(s,t)$ must exist for all given $s$ and $t$, and has transformed into a rich and well-studied subset of combinatorics over the course of the last few decades \cite{ramsey1987problem, walker1968dichromatic, McKay1995R45}. An extensive survey with current bounds, applications, and directions of the field can be found at \cite{radziszowski2011small}.

Lower bounds for a given $R(s,t,n)$ have been determined through direct mathematical proofs since the beginning of the field, and in the last 30 years saw success through the addition of computational approaches \cite{radziszowski2011small}. A counterexample graph for $R(s,t) \leq n$ is a graph of size $n$ such that there exist no cliques of size $s$ or independent sets of size $t$. $\mathcal{R}(s,t,n)$ is the set of all counterexample graphs to $R(s,t)$ for a given $n$. 

Some basic observations are as follows: Denote $k := |\mathcal{R}(s,t,n)|$. $k > 0 \Rightarrow R(s,t) > n$. $k = 0$ and $|\mathcal{R}(s,t,n-1)| > 0 \Rightarrow R(s,t) = n$. The discovery of nonempty counterexample sets has been used extensively to improve lower bounds, with $R(4,6) > 35$ and $R(5,5) > 42$ being current examples \cite{mckay1997subgraph, exoo2012ramsey}.  A full database of known counterexample graphs can be found at \cite{ramsey_graphs}. The sets of known counterexamples for $\mathcal{R}(5,5,42)$ and $\mathcal{R}(4,6,35)$ are not provably exhaustive, although McKay and Radziszowski conjecture that $\mathcal{R}(5,5,42)$ is \cite{mckay1997subgraph}. 

McKay and Radziszowski explicitly state that the known 656 $\mathcal{R}(5,5,42)$ counterexamples cannot be used to generate any $\mathcal{R}(5,5,43)$ counterexamples \cite{mckay1997subgraph}.  While no exact algorithm is provided, it can be assumed they used a variation of the one-vertex extension algorithm used in \cite{McKay1995R45}.  Exoo did not explicitly mention the use of an algorithm proving no $\mathcal{R}(4,6,36)$ could be generated from the known 37 $\mathcal{R}(4,6,35)$ counterexamples \cite{exoo2012ramsey}.

The one-vertex extension algorithm developed in \cite{McKay1995R45} has an exact description for $R(4,5)$ but would need modification for other $s$ and $t$ values. Additionally, the algorithm has a runtime bound by the number of triangles and independent sets of size 4 in each counterexample, meaning for a given counterexample the runtime is bound exponentially to $\max\{s,t\}-1$. 

\subsection{Contribution}
I provide and prove the main theorem in Section \ref{sec:main-theorem}. Using this, I first introduce a new algorithm for counterexample checking in Section \ref{sec:counter}. This algorithm has a runtime bound linearly by $\max\{s,t\}$. Then, in Section \ref{sec:ove}, I show a new algorithm for one-vertex extension, or going from $\mathcal{R}(s,t,n)$ to $\mathcal{R}(s,t,n+1)$. The algorithm has a runtime linearly bound by $\max\{s,t\}$ instead of exponentially. The algorithm also works independent of choice of $s$, $t$, and $n$. Given an input of a known $\mathcal{R}(s,t,n)$, it generates the set of all $\mathcal{R}(s,t,n+1)$ with at least 2 subgraphs within the known $\mathcal{R}(s,t,n)$, as opposed to prior algorithms that generate the set of all $\mathcal{R}(s,t,n+1)$ with at least 1 subgraph within the known $\mathcal{R}(s,t,n)$. This limits the utility of the algorithm in searching a broader space but allows for efficient execution. I validate this by verifying $\mathcal{R}(5,5,43)$ and $\mathcal{R}(4,6,37)$ are empty given current sets $\mathcal{R}(5,5,42)$ and $\mathcal{R}(4,6,36)$ on a single PC in under 15 minutes.

\section{Main Theorem}\label{sec:main-theorem}
Let $G_{n+1}$ be a simple, undirected graph of size $n+1$. 
\begin{theorem}
    If $G_{n+1}$ has $\max\{s,t\}+1$ subgraphs in $\mathcal{R}(s,t,n)$, $G_{n+1}\in \mathcal{R}(s,t,n+1)$
\end{theorem}
\begin{proof}
    Assume $G_{n+1}$ has $\max\{s,t\}+1$ subgraphs in $\mathcal{R}(s,t,n)$. Each subgraph in $\mathcal{R}(s,t,n)$ must be of size $n$, and so exclude a vertex in $G_{n+1}$. Let $S=\{v_1, \ldots, v_{\max\{s,t\}+1}\}$ denote the excluded vertices in $G_{n+1}$ for each of the $\max\{s,t\}+1$ subgraphs in $\mathcal{R}(s,t,n)$.

Assume there exists a subgraph $G'$ of $G_{n+1}$ of size $s$ or $t$. If $G'$ excludes any $v_i\in S$, then $G'$ is a subgraph of the $i$th subgraph in $\mathcal{R}(s,t,n)$. Therefore $G'$ cannot be a clique of size $s$ or an independent set of size $t$. If $G'$ includes all $v_i\in S$, then it is at least of size $\max\{s,t\}+1 > \max\{s,t\}$ which is a contradiction.
\end{proof}

\section{Counterexample Checking}\label{sec:counter}
Algorithm \ref{alg:hl-counter-checker} checks the belonging of a candidate graph $G_{n+1}$ to counterexamples $\mathcal{R}(s,t,n+1)$.

\begin{algorithm}[ht!]
    \caption{Counterexample Checking}\label{alg:hl-counter-checker}
    \begin{algorithmic}[]
        \FOR{$\max\{s,t\}+1$ subgraphs $G_n$ of size $n$ of graph $G_{n+1}$}
            \IF{$G_n$ is not isomorphic to some $G'_n\in\mathcal{R}(s,t,n)$}
                \STATE \textbf{return} false
            \ENDIF
        \ENDFOR
        \STATE \textbf{return} true
    \end{algorithmic}
\end{algorithm}

As opposed to prior algorithms exponentially related to $s$ and $t$, this algorithm requires a linear number of checks. However, other algorithms require no isomorphism checking. This algorithm requires a number of isomorphism calls linearly related to $s$ and $t$, and the isomorphism must also be found to exist or not exist from a set of $k$ possible candidates, giving an upper bound of $k(\max\{s,t\}+1)$ checks in the naive implementation.

Additionally, this algorithm requires knowledge of $\mathcal{R}(s,t,n)$. If some subset of $\mathcal{R}(s,t,n)$ is used, or the set used is not known to be exhaustive, a false return value is possible for a candidate $G_{n+1}\in \mathcal{R}(s,t,n+1)$. 

\section{One-Vertex Extension}\label{sec:ove}
Algorithm \ref{alg:hl-extension} extends from $\mathcal{R}(s,t,n)$ to $\mathcal{R}(s,t,n+1)$.

\begin{algorithm}[ht!]
    \caption{One-Vertex Extension}\label{alg:hl-extension}
    \begin{algorithmic}[0]
        \FORALL{$G_n \in \mathcal{R}(s,t,n)$}
            \STATE Define $G_{n+1} := G_n + \{v_{n+1}\}$
            \FORALL{$i\in [1,n]$}
                \STATE Graph $G_{n-1} := G_n - v_i$
                \FORALL{$G'_n \in \mathcal{R}(s,t,n)$ s.t $G'_n \neq G_n, G_{n-1} \subset G'_n$}
                    \STATE Assign edges between $G_{n-1}$ and $v_{n+1}$ s.t. $G_{n-1}+\{v_{n+1}\} = G'_n$
                    \STATE Add $G_{n+1}$ to $\mathcal{R}(s,t,n+1)$ if it is a counterexample.
                    \STATE Connect $v_{n+1}$ to $v_i$.
                    \STATE Add this new $G_{n+1}$ to $\mathcal{R}(s,t,n+1)$ if it is a counterexample.
                \ENDFOR
            \ENDFOR
        \ENDFOR
        \STATE \textbf{return} $\mathcal{R}(s,t,n+1)$
    \end{algorithmic}
\end{algorithm}

This algorithm iterates over at most $2k^2n$ candidates. Other extension algorithms look at all ways to connect a new vertex ($2^n$) to an existing counterexample, meaning they search $k2^n$ candidates. Both extension algorithms in conjunction allow for the improved choice of $\min\{2k^2n,k2^n\}$ candidates. Additionally, this algorithm explores all possible graphs with at least 2 subgraphs in $\mathcal{R}(s,t,n)$. Because $\mathcal{R}(s,t,n+1)$ is the set of all graphs with $\max\{s,t\}+1$ subgraphs in $\mathcal{R}(s,t,n)$, this exhaustively find all graphs in $\mathcal{R}(s,t,n+1)$.

Regarding checking counterexamples, existing algorithms or Algorithm \ref{alg:hl-counter-checker} can be used. When using Algorithm \ref{alg:hl-counter-checker}, only $\max\{s,t\}-1$ additional iterations are needed as two subgraphs have already been verified.

Running Algorithm \ref{alg:hl-extension} on a subset of $\mathcal{R}(s,t,n)$, only graphs with at least $\max\{s,t\}+1$ subgraphs in the subset are  generated when using Algorithm \ref{alg:hl-counter-checker} to check, and at least 2 subgraphs in the subset when checking with other schemes. This is a stricter limitation than other extension algorithms that generate graphs with at least 1 subgraph in the subset.

An additional cost of this algorithm is that, in my implementation, the condition of finding that $G_{n-1} \subset G'_n$ requires calls to an isomorphism checker.

\section{Results}
I wrote a single-threaded python implementation of the high level algorithms provided in Sections \ref{sec:counter} and \ref{sec:ove}. Code was ran on a laptop with a single intel core i7. The code, available at \url{https://github.com/aLehav/ramsey-ove-and-check}, is open-sourced and fully documented. Algorithms used and details regarding the code are in Section \ref{sec:code_imp}. I also decremented to subsets of $\mathcal{R}(s,t,n-1)$ from given subsets of some $\mathcal{R}(s,t,n)$ to broaden the search space, as done by McKay, Exoo, and Radziszowski \cite{mckay1997subgraph, exoo2012ramsey}.

\subsection{R(4,6,36)}
There are no counterexamples of size 36 with at least 7 subgraphs strictly within the 37 known $\mathcal{R}(4,6,35)$ examples \cite{ramsey_graphs}. This took 20 seconds to verify.

\subsection{R(5,5,43)}
There are no counterexamples of size 43 with at least 6 subgraphs strictly within the 656 known $\mathcal{R}(5,5,42)$ examples\cite{ramsey_graphs}. This took 9 minutes and 5 seconds to verify.

\subsection{Decrementing}
I decremented the subset of $\mathcal{R}(4,6,35)$ to a subset of $\mathcal{R}(4,6,34)$ and the subset of $\mathcal{R}(5,5,42)$ to a subset of $\mathcal{R}(5,5,41)$ and tried reconstructing to find new counterexamples. None were found.

\section{Conclusion}
I introduce a new theorem and way to look at Ramsey theory counterexamples with regard to smaller counterexamples. I provide algorithms that leverage this, and implement them to verify known results. This new approach may assist in novel methods for one-vertex extension or gluing that yield the complete sets $\mathcal{R}(4,6,36)$, $\mathcal{R}(5,5,43)$, or $\mathcal{R}(3,10,40)$ in the future.

\appendix
\section{Code Implementation}\label{sec:code_imp}
Much of the work leveraged the NetworkX library, in particular the VF2++ algorithm \cite{networkx, juttner2018vf2}. The exact libraries and data structures used may affect the runtime of these algorithms.
\subsection{Algorithms}\label{subsec:alg}
There are three bottlenecks in the practical implementation of Algorithm \ref{alg:hl-counter-checker} and Algorithm \ref{alg:hl-extension}. 
\begin{enumerate}
    \item Algorithm \ref{alg:hl-counter-checker} determines if a candidate graph has $\max\{s,t\}+1$ subgraphs of size $n$ that are isomorphic to any one of $k$ counterexamples.
    \item Algorithm \ref{alg:hl-extension} finds all $k$ counterexamples sharing a subgraph of size $n-1$ with a given candidate graph.
    \item In Algorithm \ref{alg:hl-extension}, after finding a counterexample sharing a subgraph of size $n-1$, the algorithm connects a vertex to recreate an isomorphism to this counterexample and checks if the new $G_{n+1}$ is a counterexample.
\end{enumerate}

Algorithm \ref{alg:hl-counter-checker} must determine if a candidate graph has $\max\{s,t\}+1$ subgraphs that are isomorphic to a counterexample. As mentioned, the naive implementation may iterate each of $\max\{s,t\}+1$ arbitrarily chosen subgraphs of $G_{n+1}$ over each of the $k$ counterexamples to check if they are isomorphic. However, by pre-computing the sorted edge degree of each of the $k$ counterexamples, time is saved by only checking a candidate subgraph against the subset of counterexamples with a matching sorted edge degree. This hashing significantly speeds up the code. Using the count of triangles through each vertex or number of all subgraphs of size 3 through a given vertex lowers the number of counterexamples to check against to 1 or 2 and only adds a runtime component of $O(n^3)$ per candidate.

Additionally, if for some $i\in [1,n+1]$, $G_{n+1}-\{v_i\}$ is in $\mathcal{R}(s,t,n)$, then any subgraph of $G_{n+1}$ of size $n$ must have a subgraph of size $n-1$ in $\mathcal{R}(s,t,n-1)$. This is because all subgraphs of size $n$ must have an overlapping subgraph of size $n-1$ with $G_{n+1}-\{v_i\}$. With this knowledge, we compose a mapping $\Psi$ from graphs in $\mathcal{R}(s,t,n-1)$ to sets of neighbors of the $n$th vertex that would produce a graph in $\mathcal{R}(s,t,n)$. This is done via Algorithm \ref{alg:ll-psi-gen}, described later.

Algorithm \ref{alg:hl-extension} could then work by iterating over every graph $G_n$ in $\mathcal{R}(s,t,n)$ and every index $i\in [1,n]$. For each set of neighbors in $\Psi(G_n-\{v_i\})$, $G_n - \{v_i\}$ can be connected to $v_{n+1}$ in an isomorphically equivalent way. Now, $G_{n+1}$ has two subgraphs verified to be in $\mathcal{R}(s,t,n)$. If for $\max{s,t}-1$ additional indices $j\notin \{i,n+1\}$, $\Psi(G_{n+1}-\{v_{n+1},v_j\})$ includes the neighbors of $v_{n+1}$ in $G_{n+1}-\{v_{n+1}\}$, then $G_{n+1} - \{v_j\} \in \mathcal{R}(s,t,n)$ and so $G_{n+1}$ is a counterexample.

This approach significantly reduces the number of counterexamples needed to find isomorphisms to with candidate graphs in both checking and extension. With this approach I developed the more detailed Algorithm \ref{alg:ll-extension} for one-vertex extension and Algorithm \ref{alg:ll-counter-checker} for checking within that context.

\begin{algorithm}[H]
    \caption{Detailed One-Vertex Extension}\label{alg:ll-extension}
    \begin{algorithmic}
        \STATE Generate mapping $\Psi$
        \FORALL{$G_n \in R(s,t,n)$, excluding complements when $s=t$}
            \STATE Graph $G_{n+1} := G_n + \{v_{n+1}\}$
            \FORALL{$i\in [1,n]$}
                \STATE Index $j := 1 + (i + 1 \: \text{mod}\: n+1)$
                \STATE Graph $G_{n-1} := G_n - \{v_i\}$
                \STATE Graph $G'_{n-1} := G_n - \{v_j\}$
                \STATE Generate Isomorphism $\Phi$ from $G_{n-1}$ to some graph in $\Psi$'s keys
                \STATE Generate Isomorphism $\Phi'$ from $G'_{n-1}$ to some graph in $\Psi$'s keys
                \FORALL{Sets of neighbors $M \in \Psi(\Phi(G_{n-1}))$}
                    \STATE Create edge $(v_{n+1},v_k)$ in $G_{n+1}$ for all $v_k\in \Phi^{-1}(M)$
                    \STATE $P := v_{n+1}$'s neighbors, currently equal to $\Phi^{-1}(M)$
                    \IF{$\Phi'(P\text{ excluding $v_j$})\in \Psi(\Phi'(G'_{n-1}))$}
                        \STATE Call Algorithm \ref{alg:ll-counter-checker} with all local variables
                    \ENDIF
                    \STATE Add edge $(v_{n+1}, v_i)$ to $G_{n+1}$ and $v_i$ to $P$
                    \IF{$\Phi'(P\text{ excluding $v_j$ }) \in \Psi(\Phi'(G'_{n-1}))$}
                        \STATE Call Algorithm \ref{alg:ll-counter-checker} with all local variables
                    \ENDIF
                \ENDFOR
            \ENDFOR
        \ENDFOR
        \STATE \textbf{return} all found counterexamples
    \end{algorithmic}
\end{algorithm}

As an aside, $\max_{G_{n-1}} \{|\Psi(G_{n-1}|\} \leq k$. It is usually far less than $k$, as in practice this requires all graphs in $\mathcal{R}(s,t,n)$ to share a subgraph of size $n-1$.

\begin{algorithm}[H]
    \caption{Detailed Counterexample Checking}\label{alg:ll-counter-checker}
    \begin{algorithmic}
        \STATE Tested subgraph count $t := 2$
        \STATE Index $k := j$
        \WHILE{$t \leq \max\{s,t\}$}
            \STATE $k = 1 + (k + 1 \: \text{mod}\: n+1)$
            \STATE Graph $G^*_{n-1} := G_n - \{v_k\}$
            \STATE Generate Isomorphism $\Phi^*$ from $G^*_{n-1}$ to some graph in $\Psi$'s keys
            \IF{$\Phi^*(P \text{ excluding $k$}) \notin \Psi(\Phi^*(G^*_{n-1}))$}
                \STATE $G_{n+1}$ is not a counterexample
            \ENDIF
        \ENDWHILE
        \STATE $G_{n+1}$ is a counterexample, add it to $\mathcal{R}(s,t,n+1)$
    \end{algorithmic}
\end{algorithm}

Algorithm \ref{alg:ll-psi-gen} generates $\Psi$. The reason that this algorithm must find all isomorphisms whereas the above algorithms do not is in cases where a vertex is connected via one mapping of a given counterexample but checked from another. If only one isomorphism is added to the mapping's output sets, a valid counterexample graph may incorrectly be labeled as a false negative.

\begin{algorithm}[H]
    \caption{Mapping construction}\label{alg:ll-psi-gen}
    \begin{algorithmic}
        \STATE Generate blank mapping $\Psi$
        \FORALL{$G_n \in \mathcal{R}(s,t,n)$}
            \FORALL{$i\in[1,n]$}
                \STATE $G_{n-1} := G_n - \{v_i\}$
                \STATE $M_i := v_i$'s neighbors
                \IF{$\Psi$'s domain contains a graph $G'_{n-1}$ isomorphic to $G_{n-1}$}
                    \FORALL{Isomorphisms $\Phi$ from $G_{n-1}$ to $G'_{n-1}$}
                        \STATE Add $\Phi(M_i)$ to $\Psi(\Phi(G_{n-1}))$
                    \ENDFOR
                \ELSE
                    \FORALL{Automorphisms $\Phi$ from $G_{n-1}$ to itself}
                        \STATE Add $\Phi(M_i)$ to $\Psi(\Phi(G_{n-1}))$
                    \ENDFOR
                \ENDIF
            \ENDFOR
        \ENDFOR
        \STATE \textbf{return} $\Psi$
    \end{algorithmic}
\end{algorithm}

\subsection{Runtime}
The vast majority of the runtime of this algorithm is spent on finding isomorphisms as well as determining their existence. As noted in Section \ref{subsec:alg}, hashing schemes can be used to minimize the number of counterexamples a candidate must compare against. In practice, hashing based on sorted vertex degree lowers the maximum number of valid counterexamples to 17 for $R(4,6,41)$ and 86 for $R(5,5,41)$. Hashing based on the number of triangles per vertex runs in equivalent time with NetworkX \cite{networkx} and decreases this to at most 2 valid counterexamples for both sets. Hashing based on all subgraphs of size 3 through each vertex takes measurably longer but reduces it down to 1 valid counterexample to check against for both sets. An ideal hashing algorithm would not only generate 1 valid counterexample, but in doing so also begin the isomorphic assignment. This is an important direction of future work. 

Triangle counts for a graph of size $n$ runs in $O(n^3)$. I use VF2++ \cite{juttner2018vf2} for isomorphism, making my time to hash and find the isomorphism of a single graph $h := O(n^3 + 1\cdot VF2++)$. This gives, for an arbitrary hashing scheme, $h=O(s+(C+N-1)\cdot\Phi)$ when the scheme runs in $s$ time and leaves at most $C$ valid counterexamples when finding $N$ isomorphisms from the candidate to counterexample and isomorphism generation takes $\Phi$ time for an arbitrary scheme. For dictionary construction, I aim to find the maximum $N$ per graph, and so denote it $h'$ in the following analysis.

If I extend from $\mathcal{R}(s,t,n)$ of size k, it takes $O(kh')$ to generate $\Psi$, iterate over $2kn\max_{G_{n-1}} \{|\Psi(G_{n-1}|\} = O(k^2n)$ candidate graphs, and for each take $O(\max\{s,t\}h)$ time to assign edges and check if a counterexample, for an overall runtime of $O(k^2n\max\{s,t\}h +kh')$. 

\subsection*{Acknowledgments}
Thanks to Afia Anjum for her endless patience in listening to me talk about these algorithms.

\bibliography{mybib}

\begin{thebibliography}{1}

\bibitem{ramsey1987problem}
Frank~P Ramsey.
\newblock On a problem of formal logic.
\newblock In {\em Classic Papers in Combinatorics}, pages 1--24. Springer, 1987.

\bibitem{walker1968dichromatic}
K~Walker.
\newblock Dichromatic graphs and ramsey numbers.
\newblock {\em Journal of Combinatorial Theory}, 5(3):238--243, 1968.

\bibitem{McKay1995R45}
Brendan~D. McKay and Stanislaw~P. Radziszowski.
\newblock R(4, 5) = 25.
\newblock {\em J. Graph Theory}, 19:309--322, 1995.

\bibitem{radziszowski2011small}
Stanislaw Radziszowski.
\newblock Small ramsey numbers.
\newblock {\em The electronic journal of combinatorics}, 1000:DS1--Aug, 2011.

\bibitem{mckay1997subgraph}
Brendan~D McKay and Stanis{\l}aw~P Radziszowski.
\newblock Subgraph counting identities and ramsey numbers.
\newblock {\em journal of combinatorial theory, Series B}, 69(2):193--209, 1997.

\bibitem{exoo2012ramsey}
Geoffrey Exoo.
\newblock On the ramsey number $ r (4, 6) $.
\newblock {\em the electronic journal of combinatorics}, 19(1):P66, 2012.

\bibitem{ramsey_graphs}
Brendan McKay.
\newblock Ramsey graphs, 2016.

\bibitem{networkx}
Aric Hagberg, Pieter~J Swart, and Daniel~A Schult.
\newblock Exploring network structure, dynamics, and function using networkx.
\newblock Technical report, Los Alamos National Laboratory (LANL), Los Alamos, NM (United States), 2008.

\bibitem{juttner2018vf2}
Alpár Jüttner and Péter Madarasi.
\newblock Vf2++—an improved subgraph isomorphism algorithm.
\newblock {\em Discrete Applied Mathematics}, 242:69--81, 2018.
\newblock Computational Advances in Combinatorial Optimization.

\end{thebibliography}
\bibliographystyle{unsrt}

\end{document}